\documentclass{amsart}
\usepackage{amssymb,eucal,graphics}

\theoremstyle{plain}
\newtheorem{theorem}{Theorem}
\newtheorem{lemma}{Lemma}[section]

\newtheorem{corollary}[lemma]{Corollary}
\newtheorem{examplepf}[lemma]{Example}
\newtheorem*{claim}{Claim}
\newtheorem*{stat}{\name}
\newcommand{\name}{testing}

\theoremstyle{definition}
\newtheorem{definition}[lemma]{Definition}
\newtheorem*{notation}{Notation}

\newtheorem{problem}{Problem}

\newenvironment{all}[1]{\renewcommand{\name}{#1}\begin{stat}}
                        {\end{stat}}

\newcommand{\qedc}{{\qed}~{\rm Claim.}}

\newenvironment{cproof} {\begin{proof}[Proof of Claim.]}
{\qedc\renewcommand{\qed}{}\end{proof}}

\numberwithin{equation}{section}

\newcommand{\set}[1]{\{#1\}}

\newcommand{\setm}[2]{\set{#1\mid#2}}

\newcommand{\seq}[1]{\langle#1\rangle}
\newcommand{\famm}[2]{\langle#1\mid#2\rangle}

\newcommand{\ZZ}{\mathbb{Z}}
\newcommand{\A}{\mathcal{A}}

\DeclareMathOperator{\rng}{rng}

\newcommand{\jirr}{join-ir\-re\-duc\-i\-ble}

\newcommand{\jsd}{join-sem\-i\-dis\-trib\-u\-tive}

\newcommand{\eps}{\varepsilon}
\newcommand{\vc}{\mathbin{\vee_{\mathrm{c}}}}
\newcommand{\vs}{\mathbin{\vee_{\ast}}}
\newcommand{\es}{\varnothing}

\newcommand{\ol}[1]{\overline{#1}}
\newcommand{\ootimes}{\mathbin{\bar{\otimes}}}

\newcommand{\jz}{$\langle\vee,0\rangle$}

\newcommand{\jzs}{\jz-sem\-i\-lat\-tice}

\DeclareMathOperator{\J}{J}

\newcommand{\M}{\mathbf{M}}

\newcommand{\TJ}{$(\mathrm{T}_\vee)$}

\newcommand{\FL}{\mathrm{F}_{\mathbf{L}}}

\newcommand{\bx}{\mathbf{x}}
\newcommand{\by}{\mathbf{y}}

\newcommand{\dd}{\mathrm{d}}

\begin{document}

\title[Tensor products of lattices]%
{Solutions to five problems on tensor products of lattices and
related matters}

\author[F.~Wehrung]{Friedrich Wehrung}
\address{CNRS, FRE 2271\\
         D\'epartement de Math\'ematiques\\
         Universit\'e de Caen\\
         14032 Caen Cedex\\
         France}
\email{wehrung@math.unicaen.fr}
\urladdr{http://www.math.unicaen.fr/\~{}wehrung}


\begin{abstract}
The notion of a \emph{capped tensor product}, introduced by
G.~Gr{\"a}tzer and the author, provides a convenient framework for the
study of tensor products of lattices that makes it possible to extend
many results from the finite case to the infinite case. In this paper, we
answer several open questions about tensor products of lattices. Among
the results that we obtain are the following:

\begin{all}{Theorem~\ref{T:WAmAm}}
Let $A$ be a lattice with zero. If $A\otimes L$ is a lattice for every
lattice~$L$ with zero, then $A$ is locally finite and $A\otimes L$ is a
capped tensor product for every lattice $L$ with zero.
\end{all}

\begin{all}{Theorem \ref{T:NonLCT}}
There exists an infinite, three-generated,
$2$-modular lattice $K$ with zero such that
$K\otimes K$ is a capped tensor product.
\end{all}

Here, $2$-modularity is a weaker identity than modularity, introduced
earlier by G.~Gr{\"a}tzer and the author.
\end{abstract}

\maketitle

\section{Introduction}\label{S:Intro}
For \jzs s $A$ and $B$, the tensor product $A\otimes B$ may be defined, in
a fashion formally similar to the tensor product of vector spaces in linear
algebra, as a universal object with respect to the notion of
\emph{bimorphism}, see \cite{GLQ81,GrWe2,GrWe3,GrWe4,GrWe5}.

The notion of tensor product of \jzs s becomes interesting for
\emph{lattices}, and the tensor product of two lattices is not always a
lattice. This phenomenon involves, among others, the study of
\emph{transferability} (see \cite{GGP75}), or of \emph{lower bounded
lattices} (see \cite{FJNa95}). More precisely, we say that a finite lattice
satisfies the \emph{condition \TJ}, if the relation $D_A$ of
join-dependency on the set $\J(A)$ of all \jirr\ elements of $A$
has no cycle. This condition is equivalent to saying that $A$ is a lower
bounded homomorphic image of a free lattice. For lattices $A$ and $B$ with
zero, if $A\otimes B$ is a so-called \emph{capped tensor product}, then 
$A\otimes B$ is a lattice (see Section~\ref{S:Basic}, and also
\cite{GrWe3,GrWe4}). The problem whether the converse holds is still open.
We say that a lattice $A$ with zero is \emph{amenable}, if $A\otimes L$ is
a capped tensor product, for every lattice $L$ with zero.

The following statement summarizes some of the results obtained in
\cite{GrWe3}.

\begin{theorem}\label{T:amenable}
For a lattice $A$ with zero, the following conditions are
equivalent:
\begin{enumerate}
\item $A$ is amenable.
\item $A$ is locally finite and $A\otimes B$ is a lattice, for every
lattice $B$ with zero.
\item $A$ is locally finite and $A\otimes\FL(3)$ is a lattice.
\item $A$ is locally finite and every finite sublattice of $A$ satisfies
\TJ.
\end{enumerate}
\end{theorem}

It would be nice to be able to replace amenability of $A$ by the more
straightforward condition ``$A\otimes L$ is a lattice, for every lattice
$L$ with zero'', that we shall call \emph{weak amenability} of $A$.
However, the problem whether weak amenability is equivalent
to amenability was still open at the time where Theorem~\ref{T:amenable}
was stated, as Problem~1 in \cite{GrWe3}. We solve this problem here in
the affirmative:

\begin{all}{Theorem~\ref{T:WAmAm}}
Every weakly amenable lattice is amenable.
\end{all}

In particular, the local finiteness assumption can be removed from (ii)
and (iii) in the statement of Theorem~\ref{T:amenable}.

One of the reasons why capped tensor products were introduced in
\cite{GrWe4} was to provide a wide context in which the
so-called \emph{Isomorphism Theorem} (see \cite{GrWe4}) would be valid,
thus extending the finite case established in \cite{GLQ81}. The question
whether, for lattices $A$ and $B$ with zero, $A\otimes B$ capped implies
that either $A$ or $B$ is locally finite is stated in Problem~1
in \cite{GrWe4}. We answer this question negatively, thus, at the same
time, showing the relevance of the notion of capped tensor product:

\begin{all}{Theorem \ref{T:NonLCT}}
There exists an infinite, three-generated, $2$-modular lattice $K$ such
that $K\otimes K$ is a capped tensor product.
\end{all}

The lattice $K$ of Theorem~\ref{T:NonLCT} enjoys some additional
properties. For a variety $\mathbf{V}$ of lattices, we say that a lattice
$A$ with zero is \emph{$\mathbf{V}$-amenable}, if $A\otimes L$ is a capped
tensor product, for any lattice $L$ with zero in $\mathbf{V}$. For a
positive integer $h$, let $\M^h$ denote the variety of \emph{$h$-modular
lattices}, as introduced in \cite{GrWe2}, see Section~\ref{S:Inf2Mod}; in
particular, $\M^1=\M$ is the variety of all modular lattices.

\begin{all}{Theorem \ref{T:Mamen}}
The lattice $K$ is $\M^h$-amenable for all $h>0$, although it is not
locally finite.
\end{all}

Hence the lattice $K$ provides a negative solution for Problem~5 in
\cite{GrWe3}. We also prove in Theorem~\ref{T:K2mod} that $K$ satisfies the
\emph{$2$-modular identity} introduced in \cite{GrWe2}, hence
easily solving the first half of Problem~6 in \cite{GrWe2} asking
whether the free $2$-modular lattice with three generators is finite
(\emph{answer}: no).

Finally, in Section~\ref{S:NoSimple}, we show, in particular, that
Problem~2 in \cite{GrWe3}, that asks whether there exists a simple,
nontrivial, amenable lattice has an easy negative answer.

\section{Basic concepts}\label{S:Basic}

We first recall some basic definitions about tensor products of \jzs s,
stated, for example, in \cite{GrWe4}. Let $A$ and $B$ be \jzs s. We
introduce a partial binary operation, the \emph{lateral join}, on
$A\times B$: let $\seq{a_0,b_0}$, $\seq{a_1,b_1}\in A\times B$; the
\emph{lateral join} $\seq{a_0,b_0}\vee\seq{a_1,b_1}$ is defined if
$a_0=a_1$ or $b_0 = b_1$, in which case, it is the join,
$\seq{a_0\vee a_1,b_0\vee b_1}$. A hereditary subset $I$ of $A\times B$ is
a \emph{bi-ideal} of $A\times B$, if it contains the subset
 \[
 \bot_{A,B}=(A\times\set{0_B})\cup(\set{0_A}\times B),
 \]
and it is closed under lateral joins.

The \emph{extended tensor product} of $A$ and $B$, denoted by
$A\ootimes B$, is the lattice of all bi-ideals of $A\times B$. It is easy
to see that it is an algebraic lattice. For $a \in A$ and $b \in B$, we
define $a\otimes b\in A\ootimes B$ by
 \[
 a\otimes b =\bot_{A,B}\cup
 \setm{\seq{x,y}\in A\times B}{\seq{x, y}\leqslant\seq{a, b}}
 \]
and call $a\otimes b$ a \emph{pure tensor}. A pure tensor is a
principal (that is, one-generated) bi-ideal of $A\times B$. We denote by
$A\otimes B$ the \jzs\ of all compact elements of $A\ootimes B$. It is
generated, as a \jzs, by the pure tensors. We observe that if the
semilattice $S$ of compact elements of an algebraic lattice $L$
forms a lattice in itself, then $S$ is a sublattice of $L$. In
particular, if $A\otimes B$ is a lattice, then it is a sublattice of
$A\ootimes B$.

A \emph{capping} of a bi-ideal $I$ of $A\times B$ is a subset $\Gamma$ of
$A\times B$ such that $I$ is the hereditary subset of $A\times B$
generated by $\Gamma\cup\bot_{A,B}$. We say that $I$ is \emph{capped}, if
it has a \emph{finite} capping.  A tensor product $A\otimes B$ is
\emph{capped}, if all its elements are capped bi-ideals. It is easy to see
that a capped tensor product is always a lattice.

A lattice $A$ with zero is \emph{amenable}, if $A\otimes L$ is a capped
tensor product, for every lattice $L$ with zero.

For a set $X$, we denote by $P\mapsto P^\dd$ the \emph{dualization} map on
the free lattice $\FL(X)$. We shall use several times the following result,
see Lemma~2.2(iii) of \cite{GrWe3}:

\begin{lemma}\label{L:TensUn}
Let $A$ and $B$ be lattices with zero, let $n$ be a positive integer, let
$a_0$, \dots, $a_{n-1}\in A$, $b_0$, \dots, $b_{n-1}\in B$. Then
 \[
 \bigvee_{i<n}a_i\otimes b_i=\bigcup_{P\in\FL(n)}
 P(a_0,\ldots,a_{n-1})\otimes P^\dd(b_0,\ldots,b_{n-1}).
 \]
\end{lemma}

\section{Weakly amenable lattices}

\begin{definition}\label{D:WeakAm}
A lattice $A$ with zero is \emph{weakly amenable}, if $A\otimes L$ is a
lattice, for every lattice $L$ with zero.
\end{definition}

It is obvious that every amenable lattice is weakly amenable. The question
of the converse is stated as Problem~1 in \cite{GrWe3}.
It was conjectured in \cite{GrWe3} that this problem had a negative answer.
Interestingly, this guess was too pessimistic, as we prove in
Theorem~\ref{T:WAmAm}. To prepare for this result, we first establish the
following lemma.

\begin{lemma}\label{L:PureMeet}
Let $m$, $n$ be positive integers, let $U$, $V$, $U_0$, \dots, $U_{n-1}$,
$V_0$, \dots, $V_{n-1}$ be elements of $\FL(m)$. We let
$\bx_0$, \dots, $\bx_{m-1}$, $\by_0$, \dots, $\by_{m-1}$ be the canonical
generators of the free lattice with $2m$ generators, $\FL(2m)$.
For all $R\in\FL(n)$, if the inequality
 \[
 U(\vec{\bx})\wedge V(\vec{\by})\leqslant
 R(U_j(\vec{\bx})\wedge V_j(\vec{\by})\mid j<n)
 \]
holds in $\FL(2m)$ (where we put $\vec{\bx}=\famm{\bx_i}{i<m}$ and
$\vec{\by}=\famm{\by_i}{i<m}$), then there exists a pure meet polynomial
$R^*\leqslant R$ such that
 \[
 U(\vec{\bx})\wedge V(\vec{\by})\leqslant
 R^*(U_j(\vec{\bx})\wedge V_j(\vec{\by})\mid j<n).
 \]
\end{lemma}

By a \emph{pure meet polynomial}, we mean a polynomial of the
form $\bigwedge_{i\in I}\bx_i$, where $I$ is a nonempty subset of
$\set{0,1,\dots,n-1}$.

\begin{proof}
We argue by induction on the length of $R$. If $R$ is a variable, we put
$R^*=R$. If $R=R_0\wedge R_1$, we put $R^*=R_0^*\wedge R_1^*$.

Now suppose that $R=R_0\vee R_1$, for polynomials $R_0$ and $R_1$. So the
inequality
 \[
 U(\vec{\bx})\wedge V(\vec{\by})\leqslant
 R_0(U_j(\vec{\bx})\wedge V_j(\vec{\by})\mid j<n)\vee
 R_1(U_j(\vec{\bx})\wedge V_j(\vec{\by})\mid j<n)
 \]
holds. Since the free lattice $\FL(2m)$ satisfies Whitman's condition,
one of the following inequalities holds:
 \begin{align}
 U(\vec{\bx})&\leqslant R(U_j(\vec{\bx})\wedge V_j(\vec{\by})\mid j<n)
 \label{Eq:UleqR}\\
 V(\vec{\by})&\leqslant R(U_j(\vec{\bx})\wedge V_j(\vec{\by})\mid j<n)
 \label{Eq:VleqR}\\
 U(\vec{\bx})\wedge V(\vec{\by})&\leqslant
 R_\nu(U_j(\vec{\bx})\wedge V_j(\vec{\by})\mid j<n),\qquad
 \text{for some }\nu<2.\label{Eq:UVleqR}
 \end{align}
However, \eqref{Eq:UleqR} never holds. Indeed, if $L=\FL(m)^\circ$ is the
lattice obtained by adding a new zero element (say, $0$) to $\FL(m)$, then
there exists a unique lattice homomorphism that sends $\bx_j$ to
itself and $\by_j$ to $0$ for all $j<m$, and applying that
homomorphism to \eqref{Eq:UleqR} gives the inequality $U(\vec{\bx})\leqslant 0$,
which does not hold.
Similarly, \eqref{Eq:VleqR} does not hold. So only \eqref{Eq:UVleqR}
remains, that is, there exists $\nu<2$ such that the inequality
 \[
 U(\vec{\bx})\wedge V(\vec{\by})\leqslant
 R_\nu(U_j(\vec{\bx})\wedge V_j(\vec{\by})\mid j<n)
 \]
holds. We put $R^*=R_{\nu}^*$.
\end{proof}

\begin{theorem}\label{T:WAmAm}
Every weakly amenable lattice is amenable.
\end{theorem}

\begin{proof}
Let $A$ be a weakly amenable lattice, we prove that $A$ is amenable. If
$A$ is finite, then this follows from Theorem~3 of \cite{GrWe3}.

Now the general case. Since the
class of amenable lattices with zero is closed under direct limits and
sublattices, see, for example, Theorem~2 of \cite{GrWe3}, it suffices to
prove that $A$ is \emph{locally finite}. Again by using Theorem~2 of
\cite{GrWe3}, it suffices to consider the case where $A$ is \emph{finitely
generated}, and then to prove that $A$ is finite.

Let $\vec a=\famm{a_i}{i<m}$ be a finite sequence of elements of $A$
generating $A$ as a lattice. Let $\bx_0$, \dots, $\bx_{m-1}$,
$\by_0$, \dots, $\by_{m-1}$ be the canonical generators of the free lattice
with $2m$ generators, $\FL(2m)$.

We define elements $H$ and $K$ of $A\otimes\FL(2m)$ by putting
 \begin{align*}
 H&=\bigvee_{i<m}a_i\otimes\bx_i,\\
 K&=\bigvee_{i<m}a_i\otimes\by_i.
 \end{align*}
By Lemma~\ref{L:TensUn}, the following equalities hold:
 \begin{align*}
 H&=\bigcup_{P\in\FL(m)}P(\vec a)\otimes P^\dd(\vec{\bx}),\\
 K&=\bigcup_{Q\in\FL(m)}Q(\vec a)\otimes Q^\dd(\vec{\by}),
 \end{align*}
hence
 \[
 H\cap K=\bigcup_{P,Q\in\FL(m)}(P(\vec a)\wedge Q(\vec a))\otimes
 (P^\dd(\vec{\bx})\wedge Q^\dd(\vec{\by})).
 \]
Since $A$ is weakly amenable, $A\otimes\FL(2m)$ is a lattice, hence
$H\cap K$ is a compact bi-ideal of $A\times\FL(2m)$. Thus there
are a positive integer $n$ and elements $P_0$, \dots, $P_{n-1}$,
$Q_0$, \dots, $Q_{n-1}$ of $\FL(m)$ such that the following relation
 \begin{equation}\label{Eq:PQmess}
(P(\vec a)\wedge Q(\vec a))\otimes(P^\dd(\vec{\bx})\wedge Q^\dd(\vec{\by}))
 \subseteq\bigvee_{j<n}
 (P_j(\vec a)\wedge Q_j(\vec a))
 \otimes(P_j^\dd(\vec{\bx})\wedge Q_j^\dd(\vec{\by}))
 \end{equation}
holds for all $P$, $Q\in\FL(m)$. To conclude the proof, it suffices to
prove that if $P(\vec a)$ is nonzero, then it belongs to the join closure
of $\setm{P_j(\vec a)}{j<n}$. Indeed, in that case, $|A|\leqslant 2^n$, so
$A$ is finite.

For an arbitrary $P\in\FL(m)$ such that $P(\vec a)>0$, we put $Q=P$
and we apply \eqref{Eq:PQmess}. By Lemma~\ref{L:TensUn}, there
exists $R\in\FL(n)$ such that the following system of inequalities is
satisfied:
 \begin{align}
 P(\vec a)&\leqslant R(P_j(\vec a)\wedge Q_j(\vec a)\mid j<n)
 \label{Eq:PaleRPQ}\\
 P^\dd(\vec{\bx})\wedge P^\dd(\vec{\by})&\leqslant
 R^\dd(P_j^\dd(\vec{\bx})\wedge Q_j^\dd(\vec{\by})\mid j<n).
 \label{Eq:Pdxyle}
 \end{align}
We observe that \eqref{Eq:PaleRPQ} holds in $A$, while \eqref{Eq:Pdxyle}
holds in $\FL(2m)$. By applying Lemma~\ref{L:PureMeet} to $R^\dd$ in
\eqref{Eq:Pdxyle}, we obtain a pure \emph{join} polynomial $R^*\geqslant R$
that may be substituted to $R$ in the inequality \eqref{Eq:Pdxyle}
without affecting its validity. Since $R^*\geqslant R$, the
inequality obtained by replacing $R$ by $R^*$ in
\eqref{Eq:PaleRPQ} is obviously satisfied. Therefore, we
may assume without loss of generality that \emph{$R$ is a pure join
polynomial}, that is, $R=\bigvee_{j\in J}\bx_j$ for some nonempty subset $J$
of $\set{0,1,\dots,n-1}$. By \eqref{Eq:Pdxyle}, the inequality
 \[
 P^\dd(\vec{\bx})\wedge P^\dd(\vec{\by})\leqslant
 P_j^\dd(\vec{\bx})\wedge Q_j^\dd(\vec{\by})
 \]
holds for all $j\in J$. Therefore, by substituting a new unit
element $1$ for all the $\by_j$ , we obtain that
$P^\dd\leqslant P_j^\dd$, thus
$P_j\leqslant P$. In particular, $P_j(\vec a)\leqslant P(\vec a)$
for all $j\in J$. Therefore, by \eqref{Eq:PaleRPQ}, we obtain
 \[
 P(\vec a)\leqslant\bigvee_{j\in J}(P_j(\vec a)\wedge Q_j(\vec a))\leqslant
 \bigvee_{j\in J}P_j(\vec a)\leqslant P(\vec a),
 \]
so $P(\vec a)=\bigvee_{j\in J}P_j(\vec a)$ belongs to the join closure of
$\setm{P_j(\vec a)}{j<n}$.
\end{proof}

\section{An infinite, three-generated, $2$-modular
lattice with zero}\label{S:Inf2Mod}

Let $K$ be the (infinite) lattice diagrammed on Figure~1.

\begin{figure}[htb]
\includegraphics{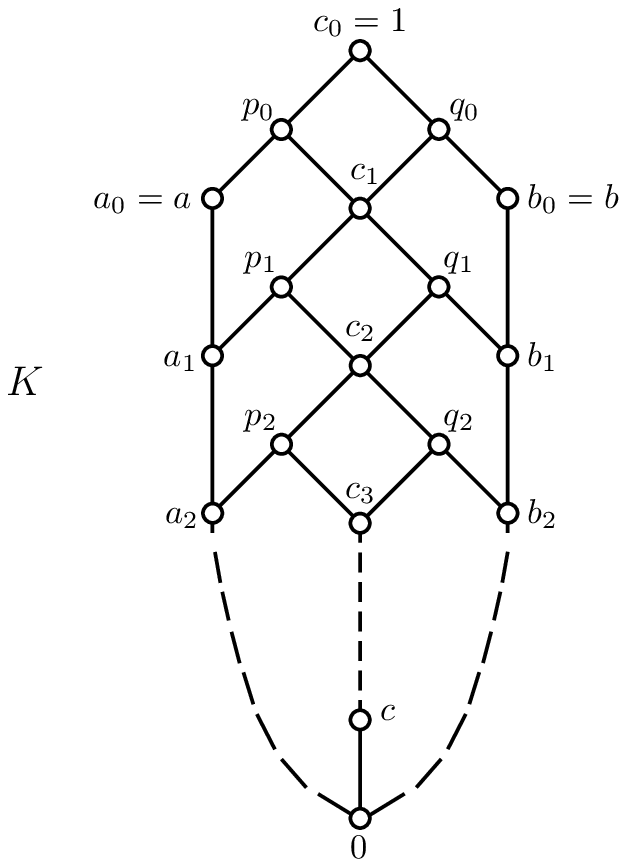}
\caption{}
\end{figure}

We observe right away the following elementary properties of $K$:

\begin{lemma}\label{L:LisJSD}
The lattice $K$ is generated by the three-element set $\set{a,b,c}$.
\end{lemma}

For any lattice $L$, we define a map $u\mapsto u^{(1)}$ from $L^3$ to
$L^3$ by the rule
 \[
 \seq{x,y,z}^{(1)}=
 \seq{x\vee(y\wedge z),y\vee(x\wedge z),z\vee(x\wedge y)},
 \quad\text{for all }x,\,y,\,z\in L.
 \]
Further, we put $u^{(0)}=u$ and $u^{(k+1)}=(u^{(k)})^{(1)}$, for all
$k<\omega$.

We say that a triple $u=\seq{x,y,z}$ of elements of $L$ is 
\begin{itemize}
\item[---] \emph{balanced}, if $u^{(1)}=u$, \emph{i.e.},
$x\wedge y=x\wedge z=y\wedge z$,

\item[---] \emph{modular}, if $\set{x,y,z}$ generates a modular
sublattice of $L$,

\item[---] \emph{distributive}, if $\set{x,y,z}$ generates a distributive
sublattice of $L$.
\end{itemize}
Of course, every distributive triple is modular.

We recall the following definition, introduced in \cite{GrWe2}:

\begin{definition}\label{D:nmod}
Let $h$ be a positive integer. A lattice $L$ is \emph{$h$-modular}, if
$u^{(h+1)}=u^{(h)}$, for any $u\in L^3$.
\end{definition}

We shall denote by $\M^h$ the variety of all $h$-modular lattices.

In particular, it is proved in \cite{GrWe2} that $1$-modularity is
equivalent to modularity. In relation to this, we recall the following
classical lemma, that says, essentially, that for every modular lattice
$M$ with zero, the tensor product $M_3\otimes M$ is capped, see
\cite{Quack}, or also \cite{Schm68} or~\cite{GrWe1}.

\begin{lemma}\label{L:u'bal}
Let $L$ be a lattice, let $u=\seq{x,y,z}$ be a modular triple of elements
of~$L$. Then $u^{(2)}=u^{(1)}$.
\end{lemma}

More generally, the following is an immediate consequence of the
definition of $h$-modularity:

\begin{lemma}\label{L:u'balh}
Let $h$ be a positive integer, let $L$ be a $h$-modular lattice, let
$u=\seq{x,y,z}$ be a triple of elements of~$L$. Then $u^{(h+1)}=u^{(h)}$.
\end{lemma}

Since the lattice $K$ contains many copies of the pentagon $N_5$ (the
five-element nonmodular lattice), it is not modular. However, it falls
relatively short of modularity:

\begin{theorem}\label{T:K2mod}
The lattice $K$ is infinite, three-generated, and
$2$-modular.
\end{theorem}

Hence, $K$ provides an answer to the second part of Problem~6 in
\cite{GrWe2}.

\begin{proof}
It remains to verify that $K$ is $2$-modular.
Let $u=\seq{x,y,z}$ be a triple of elements of $K$, we must prove that
$u^{(3)}=u^{(2)}$. If two of the elements $x$, $y$, and $z$ are
comparable, then, since $N_5$ is $2$-modular and the sublattice
generated by $\set{x,y,z}$ is a homomorphic image of $N_5$,
$u^{(3)}=u^{(2)}$ and we
are done. Hence it suffices to verify that the equality $u^{(3)}=u^{(2)}$
holds for $u$ an \emph{antichain} of $K$. If one component of $u$ is
$c$, then $u=u^{(1)}$. Otherwise, it is not hard to verify that
$u^{(1)}$ is always a triple of elements of
$K\setminus\J(K)$ (where $\J(K)$ denotes the set of all \jirr\
elements of $K$), in particular, $u^{(1)}$ is a distributive triple of
elements of~$K$. Hence, by Lemma~\ref{L:u'bal},
$(u^{(1)})^{(2)}=(u^{(1)})^{(1)}$, that is, $u^{(3)}=u^{(2)}$.
\end{proof}

\section{A non locally finite lattice that is
$\M^h$-amenable for all $h$}\label{S:MAmen}

Now let $h$ be a positive integer.
We shall prove that the lattice $K$ introduced in
Section~\ref{S:Inf2Mod} is $\M^h$-amenable, \emph{i.e.}, that $K\otimes L$
is a capped tensor product, for every $h$-modular lattice $L$ with zero.
The lattice $L$ will be fixed throughout the present section. We denote by
$\A$ the set of all maps $x\colon\J(K)\to L$ with finite range that are
\emph{antitone}, \emph{i.e.}, $p\leqslant q$ implies that $x(p)\geqslant x(q)$, for
all $p$, $q\in\J(K)$.

For any $x\in\A$, both sequences $\famm{x(a_n)}{n<\omega}$ and
$\famm{x(b_n)}{n<\omega}$ are increasing, thus, since $x$ has finite
range, the two sequences are eventually constant. We denote by
$x(a_\infty)$ and $x(b_\infty)$ their respective limits. We also denote
by $d(x)$ the least nonnegative integer that satisfies the statement
 \[
 x(a_n)=x(a_\infty)\quad\text{and}\quad x(b_n)=x(b_\infty),\qquad
 \text{for all }n<\omega\text{ such that }n\geqslant d(x).
 \]
Then we define a map $x^{(1)}$ from $\J(K)$ to $L$ by the following
equalities:
 \begin{equation}\label{Eq:Defx'}
 \begin{aligned}
 x^{(1)}(c)&=x(c)\vee(x(a_\infty)\wedge x(b_\infty)),\\
 x^{(1)}(a_0)&=x(a_0),\\
 x^{(1)}(b_0)&=x(b_0),\\
 x^{(1)}(a_{n+1})&=x(a_{n+1})\vee(x(b_n)\wedge x(c)),\\
 x^{(1)}(b_{n+1})&=x(b_{n+1})\vee(x(a_n)\wedge x(c)),\\
 \end{aligned}
 \end{equation}
for all $n<\omega$.

\begin{lemma}\label{L:x'inA}
The set $\A$ is closed under the map $x\mapsto x^{(1)}$.
\end{lemma}

\begin{proof}
For any subset $X$ of $L$, we denote by $X^\wedge$ (resp., $X^\vee$) the
meet-closure (resp., the join-closure) of $X$. Then $\rng x^{(1)}$ is a
subset of $(\rng x)^{\wedge\vee}$, hence it is finite. To conclude the
proof, it suffices to prove that $x^{(1)}(a_n)\leqslant x^{(1)}(a_{n+1})$ and
$x^{(1)}(b_n)\leqslant x^{(1)}(b_{n+1})$, for all $n<\omega$. We verify for
example the first inequality. It is trivial for $n=0$. For $n>0$, we
compute:
\begin{align*}
x^{(1)}(a_n)&=x(a_n)\vee(x(b_{n-1})\wedge x(c))\\
&\leqslant x(a_{n+1})\vee(x(b_n)\wedge x(c)) &&
(\text{because }x\text{ is antitone})\\
&=x^{(1)}(a_{n+1}),
\end{align*}
which concludes the proof.
\end{proof}

Lemma~\ref{L:x'inA} makes it possible to define inductively an
element $x^{(k)}$ of $\A$, for $x\in\A$ and $k<\omega$, by $x^{(0)}=x$
and $x^{(k+1)}=(x^{(k)})^{(1)}$, for all $k<\omega$. We further define a
map $\ell\colon\A\to L^3$ by the rule
 \[
 \ell(x)=\seq{x(a_\infty),x(b_\infty),x(c)},\qquad\text{for all }x\in\A.
 \]
We proceed with the following easy observation:

\begin{lemma}\label{L:lcommx'}
The equality $\ell(x^{(1)})=\ell(x)^{(1)}$ holds, for all $x\in\A$.
\end{lemma}

Then we put $\A'=\setm{x\in\A}{\ell(x^{(1)})=\ell(x)}$, a subset of
$\A$. We get immediately from Lemmas \ref{L:u'balh} and \ref{L:lcommx'},
together with the $h$-modularity of $L$, the following statement:

\begin{lemma}\label{L:A'closed}
The element $x^{(h)}$ belongs to $\A'$, for every $x\in\A$.
\end{lemma}

For $x\in\A'$, the formula in \eqref{Eq:Defx'} for computing
$x^{(1)}$ simplifies:
 \begin{equation}\label{Eq:Defx'onA'}
 \begin{aligned}
 x^{(1)}(c)&=x(c),\\
 x^{(1)}(a_0)&=x(a_0),\\
 x^{(1)}(b_0)&=x(b_0),\\
 x^{(1)}(a_{n+1})&=x(a_{n+1})\vee(x(b_n)\wedge x(c)),\\
 x^{(1)}(b_{n+1})&=x(b_{n+1})\vee(x(a_n)\wedge x(c)),\\
 \end{aligned}
 \end{equation}
for all $n<\omega$.
In particular, we observe that for $n\geqslant d(x)$, the inequalities
$x^{(1)}(a_n)\leqslant x(a_n)\vee(x(b_n)\wedge x(c))=
x(a_n)\vee(x(a_n)\wedge x(c))=x(a_n)$ hold, thus
$x^{(1)}(a_n)=x(a_n)=x(a_\infty)$. Similarly,
$x^{(1)}(b_n)=x(b_n)=x(b_\infty)$. Also,
$x^{(1)}(c)=x(c)$. Then an easy induction on $k$ leads to the
following result:

\begin{lemma}\label{L:d(x')=d(x)}
Let $x\in\A'$, let $n$, $k$ be nonnegative integers. If $n\geqslant d(x)$, then
$x^{(k)}(a_n)=x(a_\infty)$ and $x^{(k)}(b_n)=x(b_\infty)$. In particular,
$d(x^{(k)})\leqslant d(x)$.
\end{lemma}

The following result deals with the values of $x^{(k)}$ on the $a_n$-s and
$b_n$-s with small index $n$:

\begin{lemma}\label{L:xk+1(abn)}
Let $x\in\A'$, let $n$, $k$ be nonnegative integers. If $k\geqslant n$, then
$x^{(k+1)}(a_n)=x^{(k)}(a_n)$ and $x^{(k+1)}(b_n)=x^{(k)}(b_n)$.
\end{lemma}

\begin{proof}
We start by proving the following claim.

\begin{claim}
For any $x\in\A'$ and positive integers $n$, $k$, the following
equalities hold:
 \begin{align*}
 x^{(k)}(a_n)&=x(a_n)\vee(x^{(k-1)}(b_{n-1})\wedge x(c)),\\
 x^{(k)}(b_n)&=x(b_n)\vee(x^{(k-1)}(a_{n-1})\wedge x(c)).
 \end{align*}
\end{claim}

\begin{cproof}
We prove, for example, the first equality, by induction on $k$. It is
trivial for $k=1$. If it holds for $k$, then we compute, by using the
induction hypothesis, together with the facts that $x^{(k)}(c)=x(c)$
and $x^{(k-1)}(b_{n-1})\leqslant x^{(k)}(b_{n-1})$:
 \begin{align*}
 x^{(k+1)}(a_n)&=x^{(k)}(a_n)\vee(x^{(k)}(b_{n-1})\wedge x^{(k)}(c))\\
 &=x(a_n)\vee
 (x^{(k-1)}(b_{n-1})\wedge x(c))\vee(x^{(k)}(b_{n-1})\wedge x(c))\\
 &=x(a_n)\vee(x^{(k)}(b_{n-1})\wedge x(c)).
 \end{align*}
The proof for $x^{(k+1)}(b_n)$ is similar.
\end{cproof}

Now we prove the conclusion of Lemma~\ref{L:xk+1(abn)}, by induction on
$k$. If $k=0$, then $n=0$ and the conclusion follows from the fact that
$x^{(1)}(a_0)=x(a_0)$ and $x^{(1)}(b_0)=x(b_0)$. Now suppose that $k>0$
(and $k\geqslant n$). In the nontrivial case where $n>0$, we compute:
\begin{align*}
x^{(k+1)}(a_n)&=x(a_n)\vee(x^{(k)}(b_{n-1})\wedge x(c))&&
(\text{by the Claim above})\\
&=x(a_n)\vee(x^{(k-1)}(b_{n-1})\wedge x(c))&&
(\text{by the induction hypothesis})\\
&=x^{(k)}(a_n)&&(\text{by the Claim above}).
\end{align*}
Similarly, we could have proved that $x^{(k+1)}(b_n)=x^{(k)}(b_n)$.
\end{proof}

Now, as an immediate consequence of Lemmas \ref{L:d(x')=d(x)} and
\ref{L:xk+1(abn)}, we are able to state the following:

\begin{corollary}\label{C:xkstat}
The equality $x^{(k)}=x^{(d(x))}$ holds, for all $x\in\A'$ and all
$k<\omega$ such that $k\geqslant d(x)$.
\end{corollary}

Now we put $\A^*=\setm{x\in\A}{x^{(1)}=x}$, and $d'(x)=d(x^{(1)})+h$, for
all $x\in\A$. Hence
$\A^*\subseteq\A'\subseteq\A$. It follows from Lemmas \ref{L:A'closed}
and \ref{C:xkstat} that $x^{(k)}=x^{(d'(x))}$, for all $x\in\A$ and all
$k<\omega$ such that $k\geqslant d'(x)$. We shall denote this
element by $\tilde{x}$. Hence, $\tilde{x}$ is the least element of
$\A^*$ such that $x\leqslant\tilde{x}$, for any $x\in\A$, we shall
call it the \emph{closure} of $x$.

We denote by $\vc$ the componentwise join on $\A$, \emph{i.e.},
$(x\vc y)(p)=x(p)\vee y(p)$, for $x$, $y\in\A$ and $p\in\J(K)$. It is
clear that $\A$ is closed under $\vc$, and that it is a semilattice under
$\vc$. Hence, $\A^*$ is also a join-semilattice under componentwise
ordering, the join, that we shall denote by $\vs$, being given by
$x\vs y=\tilde{z}$ where $z=x\vc y$, for all $x$, $y\in\A^*$.

\begin{lemma}\label{L:ExtA*Emb}
Let $x$ be a map from $\J(K)$ to $L$ with finite range. Then $x$
belongs to $\A^*$ if{f} $x$ extends to a homomorphism from
$\seq{K^-,\vee}$ to $\seq{L,\wedge}$ (we put $K^-=K\setminus\set{0_K}$). Furthermore, such an extension is
unique.
\end{lemma}

\begin{proof}
For any $n<\omega$, the inequalities $c<a_n\vee b_n$, $a_{n+1}<b_n\vee c$,
and $b_{n+1}<a_n\vee c$ hold in $K$. Therefore, if $x$ extends to a
homomorphism from $K^-$ to $L$, then $x^{(1)}=x$ (see the formulas
\eqref{Eq:Defx'}).

Conversely, suppose that $x^{(1)}=x$. We prove that $x$ extends to a
unique homomorphism from $K^-$ to $L$. The uniqueness assertion is
obvious, because every element of $K^-$ is a join of finitely many, and
even at most two, elements of $\J(K)$. To prove the existence assertion,
it suffices to prove that for any elements $p$, $q$, and $r$ of $\J(K)$,
$r<p\vee q$ implies that $x(p)\wedge x(q)\leqslant x(r)$. This is obvious if
either $r\leqslant p$ or $r\leqslant q$, because $x$ is antitone. Hence suppose
that $r\nleqslant p$ and $r\nleqslant q$. We need to check the following cases:

\begin{itemize}
\item $c<a_m\vee b_n$, for $m$, $n<\omega$. Then
$x(a_m)\wedge x(b_n)\leqslant x(a_\infty)\wedge x(b_\infty)\leqslant
x^{(1)}(c)=x(c)$.

\item $a_m<b_n\vee c$, for $m$, $n<\omega$ such that $m>n$. Then
$x(b_n)\wedge x(c)\leqslant x(b_{m-1})\wedge x(c)\leqslant x^{(1)}(a_m)=x(a_m)$.

\item The case $b_m<a_n\vee c$, for $m$, $n<\omega$ such that $m>n$, is
treated similarly.
\end{itemize}
The three cases above are sufficient to conclude the proof.
\end{proof}

For an element $x$ of $\A^*$, we shall denote by $\ol{x}$ the unique
homomorphism of $K^-$ to $L$ that extends $x$. We observe that
$\rng\ol{x}=(\rng x)^\wedge$, hence $\rng\ol{x}$ is finite.

\begin{notation}
For $x\in\A$ and $I\in K\ootimes L$, let $x\nearrow I$ abbreviate the
following statement:
 \[
 \seq{p,x(p)}\in I,\quad\text{for all }p\in\J(K).
 \]
\end{notation}

\begin{lemma}\label{L:stabnea}
Let $x$, $y\in\A$, let $I\in K\ootimes L$. Then the following assertions
hold:
\begin{enumerate}
\item $x$, $y\nearrow I$ implies that $x\vc y\nearrow I$.

\item $x\nearrow I$ implies that $x^{(1)}\nearrow I$.

\item $x$, $y\in\A^*$ and $x$, $y\nearrow I$ implies that
$x\vs y\nearrow I$.
\end{enumerate}
\end{lemma}

\begin{proof}
(i) is obvious.

(ii) We assume that $x\nearrow I$, we prove that $\seq{p,x^{(1)}(p)}\in I$
for all $p\in\J(K)$. This amounts to verifying the following cases:
\begin{itemize}
\item $p=c$. From $\seq{a_n,x(a_n)}\in I$, $\seq{b_n,x(b_n)}\in I$,
$c\leqslant a_n\vee b_n$, and the fact that $I$ is a bi-ideal of $K\times L$
follows that $\seq{c,x(a_n)\wedge x(b_n)}\in I$. For $n\geqslant d(x)$, we
obtain that $\seq{c,x(a_\infty)\wedge x(b_\infty)}\in I$, hence, since
$\seq{c,x(c)}\in I$, we obtain that $\seq{c,x^{(1)}(c)}\in I$.

\item $p=a_n$, $n<\omega$. If $n=0$, then
$\seq{a_n,x^{(1)}(a_n)}=\seq{a_n,x(a_n)}\in I$. Now suppose that $n>0$.
From $\seq{b_{n-1},x(b_{n-1})}\in I$, $\seq{c,x(c)}\in I$, and
$a_n<b_{n-1}\vee c$ follows that $\seq{a_n,x(b_{n-1})\wedge x(c)}\in I$,
hence, since $\seq{a_n,x(a_n)}\in I$, we obtain that
$\seq{a_n,x^{(1)}(a_n)}\in I$.

\item $p=b_n$, $n<\omega$. This case can be treated in a similar fashion
as the previous one.
\end{itemize}

(iii) is an immediate consequence of (i) and (ii) above, together with
the fact that $x\vs y=(x\vc y)^{(n)}$ for some $n<\omega$.
\end{proof}

\begin{notation}
For $x\in\A^*$, we put
$\eps(x)=\setm{\seq{u,\xi}\in K\times L}{u>0\Rightarrow\xi\leqslant\ol{x}(u)}$.
\end{notation}

\begin{lemma}\label{L:eps}
The following assertions hold.
\begin{enumerate}
\item $\eps(x)\subseteq I$ if{f} $x\nearrow I$, for all $x\in\A^*$ and
all $I\in K\ootimes L$.

\item $\eps(x)$ is a capped element of $K\otimes L$, for any $x\in\A^*$.

\item $\eps$ is a homomorphism from $\seq{\A^*,\vs}$ to
$\seq{K\otimes L,\vee}$.
\end{enumerate}
\end{lemma}

\begin{proof}
(i) follows immediately from the fact that $\ol{x}$ is a homomorphism
from $\seq{K^-,\vee}$ to $\seq{L,\wedge}$.

(ii) It follows, again, from the fact that $\ol{x}$ is a homomorphism
from $\seq{K^-,\vee}$ to $\seq{L,\wedge}$ that $\eps(x)$ is a bi-ideal of
$K\times L$. It remains to verify that $\eps(x)$ has a finite capping. To
this end, for any $\xi\in\rng\ol{x}$, we denote by $\Gamma_\xi$ the set
of all maximal elements of $\ol{x}^{-1}\set{\xi}$. Furthermore, we put
 \[
 \Gamma=\setm{\seq{u,\xi}\in K^-\times L}
 {\xi\in\rng\ol{x}\text{ and }u\in\Gamma_\xi}.
 \]
For any $\xi\in\rng\ol{x}$, $\Gamma_\xi$ is an antichain of $K$, thus it
has at most four elements. Hence, since $\rng\ol{x}$ is finite, $\Gamma$
is finite. Now we prove that $\Gamma$ is a capping of $\eps(x)$. First,
it is obvious that $\Gamma$ is contained in $\eps(x)$. Now let
$\seq{u,\xi}\in\eps(x)$, with $u>0_K$ and $\xi>0_L$. Put
$\eta=\ol{x}(u)$. Then $u\in\ol{x}^{-1}\set{\eta}$, hence, since $K$ is
n\oe{}therian (\emph{i.e.}, every ascending chain of $K$ is eventually
constant), there exists $v\in\Gamma_\eta$ such that $u\leqslant v$. Hence,
$\seq{u,\xi}\leqslant\seq{v,\eta}$, with $\eta\in\rng\ol{x}$ and
$v\in\Gamma_\eta$, whence $\seq{v,\eta}\in\Gamma$, thus proving our
assertion. Therefore, $\eps$ maps $\A^*$ to $K\otimes L$.

(iii) It is obvious that $\eps$ is an order-preserving map from $\A^*$
(with componentwise ordering) to $K\otimes L$ (with containment).
It remains to prove that $\eps(x\vs y)\subseteq\eps(x)\vee\eps(y)$, for
all $x$, $y\in\A^*$ (the join in the right hand side is computed in
$K\otimes L$). Put $I=\eps(x)\vee\eps(y)$. Then $\eps(x)$, $\eps(y)$ are
contained in $I$, thus, by assertion (i) above, $x\nearrow I$ and
$y\nearrow I$, whence, by Lemma~\ref{L:stabnea}(iii), $x\vs y\nearrow I$,
\emph{i.e.}, by assertion (i) above, $\eps(x\vs y)\subseteq I$.
\end{proof}

To conclude the proof, we now need nothing more than a short lemma:

\begin{lemma}\label{L:puretens}
The pure tensor $u\otimes\xi$ belongs to the range of $\eps$, for all
$\seq{u,\xi}\in K\times L$.
\end{lemma}

\begin{proof}
Let $x\colon\J(K)\to L$ be the map defined by $x(p)=\xi$ if $p\leqslant u$,
$x(p)=0$ if $p\nleqslant u$, for all $p\in\J(K)$. It is easy to compute
that $\eps(x)=u\otimes\xi$.
\end{proof}

By Lemmas \ref{L:eps}(iii) and \ref{L:puretens}, the range of $\eps$
contains $K\otimes L$, while by Lemma~\ref{L:eps}(ii), every element of
the range of $\eps$ is capped. Hence, $K\otimes L$ is a capped tensor
product. Hence we have proved the following theorem:

\begin{theorem}\label{T:Mamen}
The lattice $K$ has the following properties:
\begin{enumerate}
\item $K$ is infinite, three-generated, $2$-modular.

\item $K$ is $\M^h$-amenable for all $h>0$.
\end{enumerate}
\end{theorem}

This solves Problem~5 in \cite{GrWe3} (the $2$-modularity is an
additional `luxury'). Since $K$ is $2$-modular, we obtain the
following consequence, which solves Problem~1 in~\cite{GrWe4}:

\begin{theorem}\label{T:NonLCT}
There exists an infinite, three-generated, $2$-modular lattice $K$
such that $K\otimes K$ is a capped tensor product.
\end{theorem}

\section{No simple nontrivial amenable lattices}\label{S:NoSimple}

It is proved in \cite{GrWe1} that every nontrivial lattice $L$ has a
proper congruence-preserving extension, denoted there by $M_3\langle
L\rangle$, a variant of E.T. Schmidt's $M_3[L]$ construction introduced in
\cite{Schm68}. If $L$ satisfies a certain axiom weaker than modularity,
then $M_3[L]\cong M_3\otimes L$, where $M_3$ is the modular lattice of
height two with three atoms, see
\cite{GrWe2}. The construction $M_3\otimes L$ cannot be used for general $L$
to prove that $L$ has a proper congruence-preserving extension, because it
may happen that $M_3\otimes L$ is not a lattice, see \cite{GrWe2,GrWe3}.
The basic reason for this is, of course, that $M_3$ is not amenable. This
motivated the following question:

\begin{all}{Problem~2 in \cite{GrWe3}}
Does there exist a simple, amenable lattice with more than two
elements?
\end{all}

In Proposition~9.1 of \cite{GrWe3}, we prove that no simple, amenable
(or even \jsd) lattice with more than two elements can have a largest
element. (A lattice is said to be \emph{\jsd}, if it satisfies
that $x\vee z=y\vee z$ implies that $x\vee z=(x\wedge y)\vee z$, for all
$x$, $y$, $z\in L$.) It turns out that Problem~2 in \cite{GrWe3} has a
negative answer, that follows immediately from the following easy
result:

\begin{theorem}\label{T:NoSimpleAm}
There exists no simple, locally finite lattice $S$ with more than two
elements such that any finite sublattice of $S$ has \TJ.
\end{theorem}

\begin{proof}
Suppose, towards a contradiction, that $S$ is as required. Then there are
incomparable elements $a$, $b$ of $S$. Since $S$ is simple and locally
finite, there exists a finite sublattice $L$ of $S$ such that $a$,
$b\in L$ and $\Theta_L(a\wedge b,a)=\Theta_L(a\wedge b,b)$
($\Theta_L(x,y)$ denotes the principal congruence of $L$ generated by
the pair $\seq{x,y}$). By assumption, $L$ satisfies \TJ. Hence, $L$
satisfies the statement, denoted in \cite{AdGo} by (DPT), that
 \[
 \Theta(u_0,u)=\Theta(v_0,v)\quad\text{implies that}\quad
 u\wedge v\nleqslant u_0\text{ and }u\wedge v\nleqslant v_0,
 \]
for all $u_0<u$ and $v_0<v$ in $L$, see \cite[p.~73]{Day}.
Putting $u=a$, $v=b$, and $u_0=v_0=a\wedge b$, we obtain a
contradiction.
\end{proof}

In contrast with Theorem~\ref{T:NoSimpleAm}, we observe the following
example:

\begin{examplepf}\label{Ex:SJSD}
There exists an infinite, simple, locally finite,
\jsd\ lattice with zero.
\end{examplepf}

\begin{proof}
Consider the lattice $S$ of all bounded intervals of the chain $\ZZ$ of all
integers, partially ordered under containment. Then it is well-known (and
easy to verify directly) that $S$ is locally finite and
\jsd. Since $S$ is atomistic (that is, every element of
$S$ is a join of finitely many---in fact, two---atoms), in order to prove
that $S$ is simple, it suffices to prove that
$\Theta_S(\es,a)=\Theta_S(\es,b)$, for any atoms $a$ and $b$ of $S$ such
that $a\ne b$.

Observe that the atoms of $L$ are exactly the
singletons of the form $\set{n}$, for $n\in\ZZ$. So there are $u$,
$v\in\ZZ$ such that $a=\set{u}$ and
$b=\set{v}$. Without loss of generality, we may assume
that $u<v$. Pick $x$, $y\in\ZZ$ such that $x<u<v<y$. Then
$\set{u}<\set{v}\vee\set{x}$ and
$\set{v}<\set{u}\vee\set{y}$. Since $\set{u}$, $\set{v}$, $\set{x}$, and
$\set{y}$ are distinct atoms of $S$, it follows that
$\Theta_S(\es,\set{u})=\Theta_S(\es,\set{v})$. Therefore, $S$ is simple.
\end{proof}

We observe a difference between these results and the
easy observation that states that there is no nontrivial simple,
\jsd\ lattice with a largest element, see Proposition~9.1
of \cite{GrWe3}. Namely, the proof of Theorem~\ref{T:NoSimpleAm} requires
amenability, which is necessary in view of
Example~\ref{Ex:SJSD}, while Proposition~9.1 of
\cite{GrWe3} requires only join-semidistributivity.

\section{New open problems}\label{S:Pbs}

\begin{problem}\label{Pb:VarWAmAm}
Let $\mathbf{V}$ be a variety of lattices, let $A$ be a lattice with zero.
If $A\otimes L$ is a lattice for any $L\in\mathbf{V}$, is $A$
$\mathbf{V}$-amenable?
\end{problem}

Solving Problem~\ref{Pb:VarWAmAm}, even for a given variety $\mathbf{V}$
(for example, the variety $\M$ of all modular lattices), may also provide
some insight towards a solution of the (still open) Problem~3 in
\cite{GrWe4}, that asks whether every tensor product of lattices that is a
lattice is capped.

By Theorem~\ref{T:WAmAm}, Problem~\ref{Pb:VarWAmAm} has a positive
solution for $\mathbf{V}=\mathbf{L}$, the variety of all lattices.

\begin{problem}
Let $\mathbf{V}$ be a nontrivial variety of lattices. Does there exist a
non locally finite, $\mathbf{V}$-amenable lattice?
\end{problem}

In Theorem~\ref{T:Mamen}, we prove that there exists a non locally
finite lattice that is $\M^h$-amenable for all $h>0$.

\section*{Acknowledgment}

The author is grateful to Marina Semenova for having read the paper and
pointed several oversights and misprints, and to the referee for his
helpful comments, in particular the observation that $K\otimes K$ is a
capped tensor product.


\begin{thebibliography}{99}

\bibitem{AdGo}
K.V. Adaricheva and V.A. Gorbunov,
\emph{On lower bounded lattices},
Algebra Universalis \textbf{46} (2001), 203--213.

\bibitem{Day}
A. Day,
\emph{Characterizations of finite lattices that are bounded homomorphic
images or sublattices of free lattices},
Canad. J. Math. \textbf{31} (1979), 69--78.

\bibitem{FJNa95}
R. Freese, J. Je\v{z}ek, and J.B. Nation,
\emph{Free Lattices},
Mathematical Surveys and Monographs, Vol.
\textbf{42}. American Mathematical Society, Providence, RI, 1995.
viii+293 pp.

\bibitem{GGP75}
H.S. Gaskill, G. Gr\"atzer, and C.R. Platt,
\emph{Sharply transferable lattices},
Canad.\ J. Math. \textbf{27} (1975), 1246--1262.

\bibitem{Grat98}
G. Gr{\"a}tzer,
\emph{General Lattice Theory. Second Edition},
Birkh\"auser Verlag, Basel. 1998. xix+663~pp.

\bibitem{GLQ81}
G. Gr\"atzer, H. Lakser, and R.W. Quackenbush,
\emph{The structure of tensor products of semilattices with zero},
Trans. Amer. Math. Soc. \textbf{267} (1981), 503--515.

\bibitem{GrWe1}
G. Gr\"atzer and F. Wehrung,
\emph{Proper congruence-preserving extensions of lattices},
Acta Math. Hungar. \textbf{85} (1999), 169--179.

\bibitem{GrWe2}
\bysame,
\emph{The $M_3[D]$ construction and $n$-modularity},
Algebra Universalis \textbf{41}, no. 2 (1999), 87--114.

\bibitem{GrWe3}
\bysame,
\emph{Tensor products and transferability of semilattices},
Canad. J. Math. \textbf{51}, no.~4 (1999), 792--815.

\bibitem{GrWe4}
\bysame,
\emph{Tensor products of semilattices with zero, revisited},
J. Pure Appl. Algebra \textbf{147} (2000), 273--301.

\bibitem{GrWe5}
\bysame,
\emph{A survey of tensor products and related constructions in two
lectures}, Algebra Universalis \textbf{45} (2001), 117--134.

\bibitem{Quack}
R.W. Quackenbush,
\emph{Non-modular varieties of semimodular lattices with spanning
$M_3$}, Discr. Math. \textbf{53} (1985), 193--205.

\bibitem{Schm68}
E.T. Schmidt,
 \emph{Zur Charakterisierung der Kongruenzverb\"{a}nde der Verb\"{a}nde},
 Mat. \v{C}asopis Sloven. Akad. Vied \textbf{18} (1968), 3--20.

\end{thebibliography}
\end{document}